\documentclass{cedi}
\pdfoutput=1
\usepackage[T1]{fontenc}
\usepackage{graphicx}
\usepackage{flushend}
\usepackage{amssymb}
\usepackage{amsthm}

\usepackage{alltt}
\usepackage{graphicx}
\usepackage{enumerate}

\newcommand{\labelof}[0]{\lambda}

\newcommand{\adjvertices}[0]{N}
\newcommand{\adjedges}[0]{N_{e}}
\newcommand{\adjfaces}[0]{N_{f}}
\newcommand{\weight}[0]{w}

\newcommand{\Zset}[0]{\mathbb{Z}}
\newcommand{\Nset}[0]{\mathbb{N}}
\newcommand{\genand}[1]{\left \{ \begin{array}{c} \begin{array}{l}\displaystyle #1 \end{array} \end{array} \right .}

\newenvironment{algorithm}[1]
 {\centerline{\small \stepcounter{algoritmo}\textbf{Algorithm \arabic{algoritmo}: }#1} \smallskip} 

\begin{document}

\title{\textbf{A Heuristic for Magic and Antimagic Graph Labellings}}

\normalsize

\renewcommand{\labelitemi}{\textbullet}

\author{\hspace*{-0.5cm}
\begin{tabular}{@{\extracolsep{0mm}}c@{\extracolsep{2mm}}c@{\extracolsep{2mm}}c}
Fran\c{c}ois Bertault \thanks{Work done while visiting the University of Newcastle, Australia} & Ramiro Feria-Pur\'on \\
\small{Tom Sawyer Software} \\
\small{181 Montecito Avenue} & \small{Exper. Station "Indio Hatuey"} \\
\small{Oakland, CA 95610, USA} & \small{Univ. of Matanzas, CP 44280, Cuba} \\
\small{bertault@tomsawyer.com} & \small{ramiro.feria.puron@gmail.com} 
\\\\
Mirka Miller \thanks{Also University of West Bohemia, Czech Republic, and Kings College, London, UK} & Elaheh Vaezpour \thanks{Work done while visiting the University of Newcastle, Australia} \\
Hebert P\'erez-Ros\'es  \thanks{Work done while visiting the University Jaume I, Castellon, Spain} \\
\small{School of Electrical Eng. and Comp. Science} & \small{Dept. of Computer Eng.} \\
\small{The University of Newcastle} & \small{Amirkabir Univ. of Technology} \\
\small{Callaghan 2308 NSW, Australia} & \small{Tehran, Iran} \\
\small{ \{Mirka.Miller, Hebert.Perez-Roses\}@newcastle.edu.au} & \small{e\_vaezpour@yahoo.com} \\
\end{tabular}
}

\date{}

\maketitle

\thispagestyle{empty}

\begin{abstract}
\noindent Graph labellings have been a very fruitful area of research in the last four decades. However, despite the staggering number of papers published in the field (over 1000), few general results are available, and most papers deal with particular classes of graphs and methods. Here we approach the problem from the computational viewpoint, and in a quite general way. We present the existence problem of a particular labelling as a combinatorial optimization problem, then we discuss the possible strategies to solve it, and finally we present a heuristic for finding different classes of labellings, like vertex-, edge-, or face-magic, and $(a, d)$-antimagic $(v, e, f)$-labellings. The algorithm has been implemented in C++ and MATLAB, and with its aid we have been able to derive new results for some classes of graphs, in particular, vertex-antimagic edge labellings for small graphs of the type $P_2^r \times P_3^s$, for which no general construction is known so far. 
\end{abstract}


\section{Introduction}
\label{sec:intro}

Interest in graph labelling problems began in the mid-1960s, and in the subsequent four decades well over 1000 papers on this topic have appeared. Labelled graphs serve as useful models for a broad range of applications such as coding theory, x-ray crystallography, radar, astronomy, circuit design, and communication networks \cite{BloGo77}. More recent applications to frequency allocation and image authentication are described in \cite{Fiala08} and \cite{Chan09}, respectively. The central problem here lies in finding new classes of finite (or infinite) graphs labelled in certain ways, or in constructing new types of labellings, under the assumption that it is very difficult to develop a closed theory giving a complete characterization of the set of graphs permitting a new labelling. For a thorough treatment of the different types of labelling, latest results, and many open problems concerning them, see \cite{Baca-Miller-monograph}, \cite{Gallian-survey}, and \cite{Wallis-book}. 
\\
In spite of the wide literature on labellings, there are relatively few general results, since most papers deal with constructive methods that are only applicable to a restricted class of graphs. Moreover, very few authors have adressed the problem from the computational viewpoint. Notable exceptions are \cite{McKay98}, \cite{Ex02}, \cite{Fiala05}, \cite{Fiala08}, \cite{KK96}, \cite{Sun94}. In particular, \cite{Sun94} presents a general algorithm to determine whether a given graph has a magic labelling, and constructs the labelling accordingly. The algorithm consists of an exhaustive backtrack search with different mechanisms to prune the search space. However, no complexity analysis is presented in the paper, and there is no implementation of the algorithm as to date.  
\\
In this article we also address the labelling problem from the computational viewpoint, and we develop a very general algorithm for implementing various graph labelling methods. Our approach consists of reducing a combinatorial labelling problem to an optimization problem, where the objective function depends on the type of labelling that we seek. Thus, a very general algorithm for different types of labelling can be formulated. In particular, our algorithm can be used for magic, antimagic, and $(a,d)$-antimagic labellings. 
\\
We start by defining the basic concepts and stating our terminology (Section \ref{sec:definitions}), then we discuss different strategies for solving the associated combinatorial optimization problem (Section \ref{sec:strategies}), and finally we develop one of these strategies in more detail, and discuss some experimental results obtained with that particular strategy (Section \ref{sec:experimental}). 

\section{Terminology and definitions}
\label{sec:definitions}

A \emph{graph} $G=(V, E, F)$ is defined by a finite set $V$ of \emph{vertices}, a finite set $E$ of \emph{edges} and a finite set $F$ of \emph{faces}. An edge is an unordered pair of vertices of $V$. If $e=\{u,v\}\in E$, $u$ and $v$ are  \emph{adjacent} and $e$ is \emph{incident} to $u$ and $v$. A face $f$ is a sequence of vertices $f=v_0, \dots, v_{k}$ such that $v_0 = v_k$ and $\forall i = 0, \dots, k-1, \{ v_i, v_{i+1} \} \in E$. The face $f$ is said to be \emph{incident} to each vertex $v_i, i=1,\dots, k$, and also incident to the edges $e = \{ v_i, v_{i+1} \}, i=0, \dots, k-1$. Two faces are said to be \emph{adjacent} if they are incident to a same edge. For a vertex $v$ (resp. an edge $e$, resp. a face $f$), we denote by $\adjvertices(v)$ the set of vertices adjacent to $v$ (resp. by $\adjvertices(e)$ the set of vertices that $e$ is incident to, resp. by $\adjvertices(f)$ the set of vertices that $f$ is incident to), by $\adjedges(v)$ the set of edges incident to $v$ (resp. by $\adjedges(f)$ the set of edges that $f$ is incident to), and by $\adjfaces(v)$ the set of faces incident to $v$ (resp. by $\adjfaces(e)$ the set of faces incident to $e$, resp. $\adjfaces(f)$ the set of faces adjacent to $f$).

The \emph{weight} $\weight(v)$ of a vertex $v\in V$ under a \emph{labelling function} $\labelof : V \cup E \cup F \rightarrow \Zset$ is defined by:
\begin{equation}\weight(v)=\labelof(v)+ \sum_{e \in \adjedges(v)} \labelof(e) +\sum_{f\in \adjfaces(v)} \labelof(f)
\end{equation}

The weight $\weight(e)$ of an edge $e=\{u, v\}\in E$ under $\labelof : V \cup E \cup F \rightarrow \Zset$ is defined by:
\begin{equation}
\weight(e)=\labelof(e)+\labelof(u)+\labelof(v) +\sum_{f \in \adjfaces(e)} \labelof(f)
\end{equation}

The weight $\weight(f)$ of a face $f\in F$ under $\labelof : V \cup E \cup F \rightarrow \Zset$ is defined by:
\begin{equation}
\weight(f)=\labelof(f)+  \sum_{v \in \adjvertices(f)} \labelof(v) + \sum_{e\in \adjedges(f)} \labelof(e) 
\end{equation}

For any set $S$, we define $\delta_{0}(S)=\emptyset$ and $\delta_{1}(S)=S$. Then $\labelof$ is a \emph{$(v,e,f)$-labelling function}, with $v \in \{ 0,1 \}, e \in \{ 0,1 \}$ and $f \in \{ 0,1 \}$ if $\labelof$ is a bijection $\labelof : \delta_{v}(V) \cup \delta_{e}(E)\cup \delta_{f}(F)  \rightarrow \{1,2,\dots,|\delta_{v}(V) \cup \delta_{e}(E)\cup \delta_{f}(F)|\}$. 

A $(1,0,0)$-labelling is called a \emph{vertex-labelling}, a $(0,1,0)$-labelling an \emph{edge labelling}, a $(0,0,1)$-labelling a \emph{face labelling}, and a $(1,1,0)$-labelling a \emph{total labelling}. A $(1,e,f)$-labelling with $\labelof(v)\in[1,|V|]$ for every vertex $v$ is a \emph{super labelling}. 

If there exists $k \in \Nset$ such that $\forall u\in V$ (resp. $\forall u\in E$ resp. $\forall u\in F$) we have  $\weight(u)=k$ under a labelling function $\labelof$, then $\labelof$ is called a \emph{vertex-magic} (resp. \emph{edge-magic}, resp. \emph{face-magic}) labelling. The integer $k$ is called the \emph{magic constant} of the labelling. 

If, on the contrary, $\forall u,v\in V$ (resp. $\forall u,v\in E$, resp. $\forall u,v\in F$) we have  $\weight(u) \neq \weight(v)$ under a labelling function $\labelof$, then $\labelof$ is called a \emph{vertex-antimagic} (resp. \emph{edge-antimagic}, resp. \emph{face-antimagic}) labelling.  

A particular case of antimagic labellings occurs if there exists $(a,d) \in \Nset \times \Nset$ such that $\forall i=0,\dots, |V|-1$ (resp. $\forall i=0,\dots, |E|-1$, resp. $\forall i=0,\dots, |F|-1$)  there is a unique $u \in V$ (resp. unique $u \in E$, resp. unique $u \in F$) such that $\weight(u)=a+i \times d$ under the labelling function $\labelof$. Then $\labelof$ is called a \emph{vertex-$(a,d)$-antimagic} (resp. \emph{edge-$(a,d)$-antimagic}, resp. \emph{face-$(a,d)$-antimagic}) labelling. $(a,d)$-antimagic labellings were first introduced in \cite{vertex-antimagic}. 

\begin{figure}[htbp]
\begin{center}
	 	\includegraphics[width=0.4\textwidth]{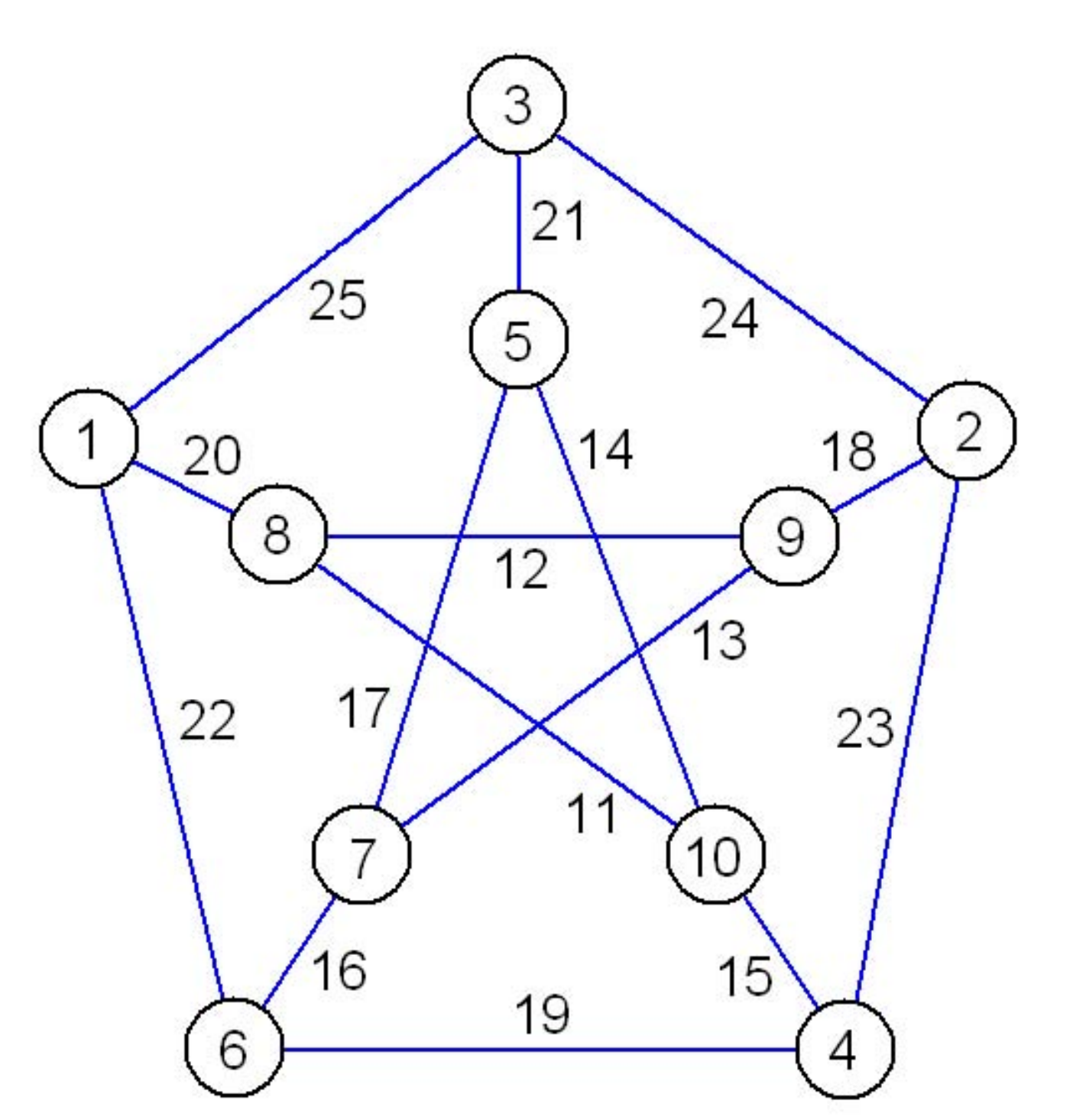}
	 \caption{Super edge-magic total labelling with magic constant 29 of the Petersen graph $P(5,2)$.}
	 \label{fig1}
\end{center}	
\end{figure}

Figure \ref{fig1} represents a super edge-magic total labelling of the Petersen graph.

By a graph labelling problem we mean the problem of assigning labels to vertices, edges and/or faces of a graph in such a way that the properties of a specified labelling are satisfied.

\section{Different strategies for finding labellings}
\label{sec:strategies}

The basic idea behind our approach is to view the problem of finding a labelling as a global optimization problem. Given an optimization function $\phi: \Nset^n \rightarrow \Nset$, the aim is to find the value $x$ for which $\phi$ is minimized, under the following constraints:
\begin{eqnarray}
\label{cpbm} 
\genand{x=[x_1,\dots,x_n]^t \\
1 \leq x_i \leq n \\
\forall i \neq j, x_i \neq x_j}
\end{eqnarray}

Several global optimization methods are available to solve discrete optimization problems \cite{BlumRoli03}. We discuss three possible solutions for our combinatorial optimization problem. The first one is based on an integer linear programming formulation, combined with the branch and bound technique \cite{Ch83}. The second one is a local search approach (more precisely, a variant of simulated annealing) \cite{PrFlTeVe86}. Finally, we discuss evolutionary techniques, as a way to achieve robustness in combination with local search.  

\subsection{Integer linear programming}
\label{sub:ilp}

Let $U = \delta \sb{v}(V) \cup \delta\sb{e}(E) \cup \delta\sb{f}(F)$, and $n=|U|$. To each element $u$ of $U$ we associate a distinct integer $\beta(u)$, in the interval $[1,n]$, called the {\it index} of $u$. Then we define $n^2$ variables $x_{ij}, i=1,\dots,n, j=1,\dots, n$. Here $x_{ij}=1$ means that the label $\labelof(u)$ of the element $u$ with index $i$ is $j$. Since an element of the set $U$ has a unique label, we have the following constraints:
\begin{eqnarray}
\genand{\forall i,j \in [1,\dots,n], x_{ij} \in \{ 0,1 \} \\
\forall i, \sum_{j=1,\dots,n} x_{ij}=1}
\end{eqnarray}

To ensure that every label is distinct, we add the following constraint:

\begin{eqnarray}
\sum_{j=1,\dots,n} (n+1)^{j-1} \sum_{i=1,\dots,n} x_{ij} = \nonumber \\
= \sum_{i=j,\dots,n} (n+1)^{n-1}
\end{eqnarray} 

The optimization function depends on the type of labelling considered. For example, for an edge-magic total labelling with constant sum $K$, the optimization function is defined by:
\begin{equation}
\phi(x)=\sum_{e\in E} |\phi(e,x)| \label{ofabs}
\end{equation}

where, for an edge $e = \{ u, v \}$

\begin{eqnarray}
\phi(e,x)=\sum_{j=1,\dots,n} j x_{\beta(u)j} + \sum_{j=1,\dots,n} j x_{\beta(v)j} + \nonumber \\ 
+ \sum_{j=1,\dots,n} j x_{\beta(e)j}-K 
\end{eqnarray} 

The optimization function $\phi$ contains an absolute value. An equivalent linear programming problem is obtained by adding a new variable $y_i$ for each edge $e_i = \{ u_i, v_i \}$, and by adding the following constraints:

\begin{eqnarray}
\genand{y_i + \phi(e_i,x) \geq 0 \\
y_i - \phi(e_i,x)\geq 0 \\ }
\end{eqnarray}

Minimizing the function $\phi$ defined in Equation \ref{ofabs} is equivalent to minimizing 

\begin{equation}
g(x)=\sum_{1=1,\dots,|E|} y_i
\end{equation}

The main problem of this method is numerical instability, and the quadratic number of variables and constraints relative to the number of labels. Therefore, this method is limited in practice to small graphs.

\subsection{Heuristic based on simulated annealing}
\label{sub:annealing}

Our second solution to the combinatorial optimization problem (\ref{cpbm}) is a local search algorithm; more precisely, a simplified variant of simulated annealing. The standard simulated annealing method is modified so as to preserve the constraint of distinctness of labels, and so that only simple arithmetic operations are required at each step. 
\\
The method proceeds as follows: We start with an initial configuration of labels that satisfies the constraints. At each step, we perform a local modification of the current configuration, and we evaluate the objective function for that new configuration. If the value is smaller than the previous one, we accept the new configuration. If the value is greater, we accept the new configuration with probability $q$, only if we failed to find a better solution in the $p$ previous attempts.
\\
There are several local changes that can be made to a particular labelling, provided that the distinctness of the labels is preserved. The simplest of these transformations is interchanging two labels, but more sophisticated transformations can be devised, such as any permutation of the labels incident to a given face, or any permutation of all the labels. The advantage of the exchange rule is not only its simplicity, but also that in principle, all possible labellings can be explored by repeated exchanges, since the symmetric group of all permutations is generated by transpositions.  
\\
The objective function $\phi$ to be minimized in our combinatorial problem depends on the type of labelling that we are looking for. The function estimates, given labels on nodes, edges and faces, satisfying the constraints, if the weights of the nodes, edges or faces of the graph are close to a solution of the magic or antimagic labelling problem. 
\\
In the case of a magic labelling problem, the objective function is very simple: it corresponds to the standard deviation of the weights. More precisely, the $S$-magic optimization function $f_S$ of a set $S \in \{ V,E,F \}$, under a labelling $\labelof$, is defined as 

\begin{equation}
f_S(G)=\sum_{s\in S} (\weight(s)-\bar{\weight}_S(G))^2
\end{equation}

\begin{equation}
\mbox{ where } \bar{\weight}_S(G)=\lceil \frac{1}{|S|} \sum_{s \in S} \weight(s)  \rceil 
\end{equation}

Note that $f_V(G)=0$ (resp. $f_E(G)=0$, resp. $f_F(G)=0$) if, and only if, $\labelof$ is a vertex-magic (resp. edge-magic, resp. face-magic) labelling. In each case, $\bar{\weight}_S(G)$ corresponds to the magic constant of the labelling.

The antimagic optimisation function $g_S$ of a  set $S \in \{ V,E,F \}$ under a labelling $\labelof$ is defined as follows. We order the elements in $S$ according to their weight: $s_1 \leq s_2 \leq \dots \leq s_{|S|}$, with $s_i\leq s_{i+1}$ if $\weight(s_i)\leq\weight(s_{i+1})$, and then we simply count the number of duplicates (i.e., the number of $i$'s for which $s_i = s_{i+1}$). Obviously, $g_S(G)=0$ iff $\labelof$ is an antimagic labelling of $G$ (be it a vertex-, edge-, or face-antimagic labelling). 

The above \lq duplicates' function, however, cannot be used for the case of $S-(a,d)$-antimagic labellings. For the $S-(a,d)$-antimagic optimisation function $h_S$ of $S\in \{ V,E,F \}$ under $\labelof$, we order the elements of $S$ according to their weight, as before, and then

\begin{eqnarray}
h_S(G) = (\weight(s_1) - a)^2 + \nonumber \\
+ \sum_{ i = 1,\ldots, |S|} ( \weight(s_{i+1}) - \weight(s_i) - d)^2
\end{eqnarray}

Again, note that $h_V(G)=0$ (resp. $h_E(G)=0$, resp. $h_F(G)=0$) when $\labelof$ is a vertex-$(a,d)$-antimagic (resp. edge-$(a,d)$-antimagic, resp. face-$(a,d)$-antimagic) labelling.

\section{Implementation and experimental results}
\label{sec:experimental}

Algorithm 1 formalizes the ideas discussed in \ref{sub:annealing}. We have chosen the simulated annealing approach mainly because of its simplicity. The algorithm receives as input the graph $G$, the objective function $\phi$ to be minimized, and some additional control parameters. The algorithm searches for a minimum value of $\phi$ (hopefully the global minimum, which is 0), and returns the corresponding labelling that yields this minimum. If the final value of $\phi$ is 0, then we get the desired labelling (magic, antimagic, etc.). The parameters $p$ and $q$ influence the number of iterations and the average time needed to find a solution, but they usually do not require fine tuning. Reasonable choices for these parameters are  $p=m(m-1)/2$ (the number of edge transpositions), and $q = 2/p$. 
\\
==================== \\
\noindent 
\texttt{
\noindent \underline{Input:} A graph $G=(V,E,F)$, \\
three booleans: $v$, $e$ and $f$,  \\
an objective function \\
$\phi \in \{ f_V, f_E, f_F, g_V, g_E, g_F, h_V, h_E, h_F \}$, \\
a positive integer $p$, \\
and a real value $0 < q < 1$. \\\\
\underline{Output:} If $\phi(G)=0$, then $\labelof$ is a \\
vertex, edge or face-magic or antimagic \\
$(v, e, f)$-labelling, depending on the \\
choice of the function $\phi$. 
\begin{tabbing} {xxxx} \= {xxxx} \= {xxxx} \= {xxxx} \kill
$U$ := $\delta_v(V) \cup \delta_e(E) \cup \delta_f(F)$; \\
$n$ := $|U|$; \\
Assign a distinct label $\labelof(u)$ in the range \\
$[1,n]$ to each element $u$ of $U$; \\
\underline{For all} $u \in V \cup E \cup F - U $ \underline{do} $\labelof(u)$:= $0$; \\
val := $\phi(G)$; \\
\underline{While} val $\neq 0$ \underline{do} \\
  \> Choose randomly two elements $r, s \in U$; \\
  \> $\labelof(r)$ :=: $\labelof(s)$; \\
  \> new\_val := $\phi(G)$; \\
  \> \underline{If} new\_val $<$ val \underline{then} \\
  \> \>  nb\_miss := 0; \\
  \> \>  val := new\_val; \\
  \> \underline{else if} (nb\_miss $> p$) and \\
  \> \> \> (random(0..1) $\leq q$) \underline{then} \\
  \> \>  nb\_miss := 0; \\
  \> \>  val := new\_val; \\  
  \> \underline{else} \\
  \> \> nb\_miss := nb\_miss + 1; \\
  \> \> $\labelof(r)$ :=: $\labelof(s)$; \\  
  \> \underline{end} \\
\underline{end}  
\end{tabbing}}
\begin{algorithm}
{Simulated annealing labelling \\}
\label{simulatedannealing}  
\end{algorithm}
==================== \\
 
\begin{figure}[htbp]
\begin{center}
	 	\includegraphics[width=0.5\textwidth]{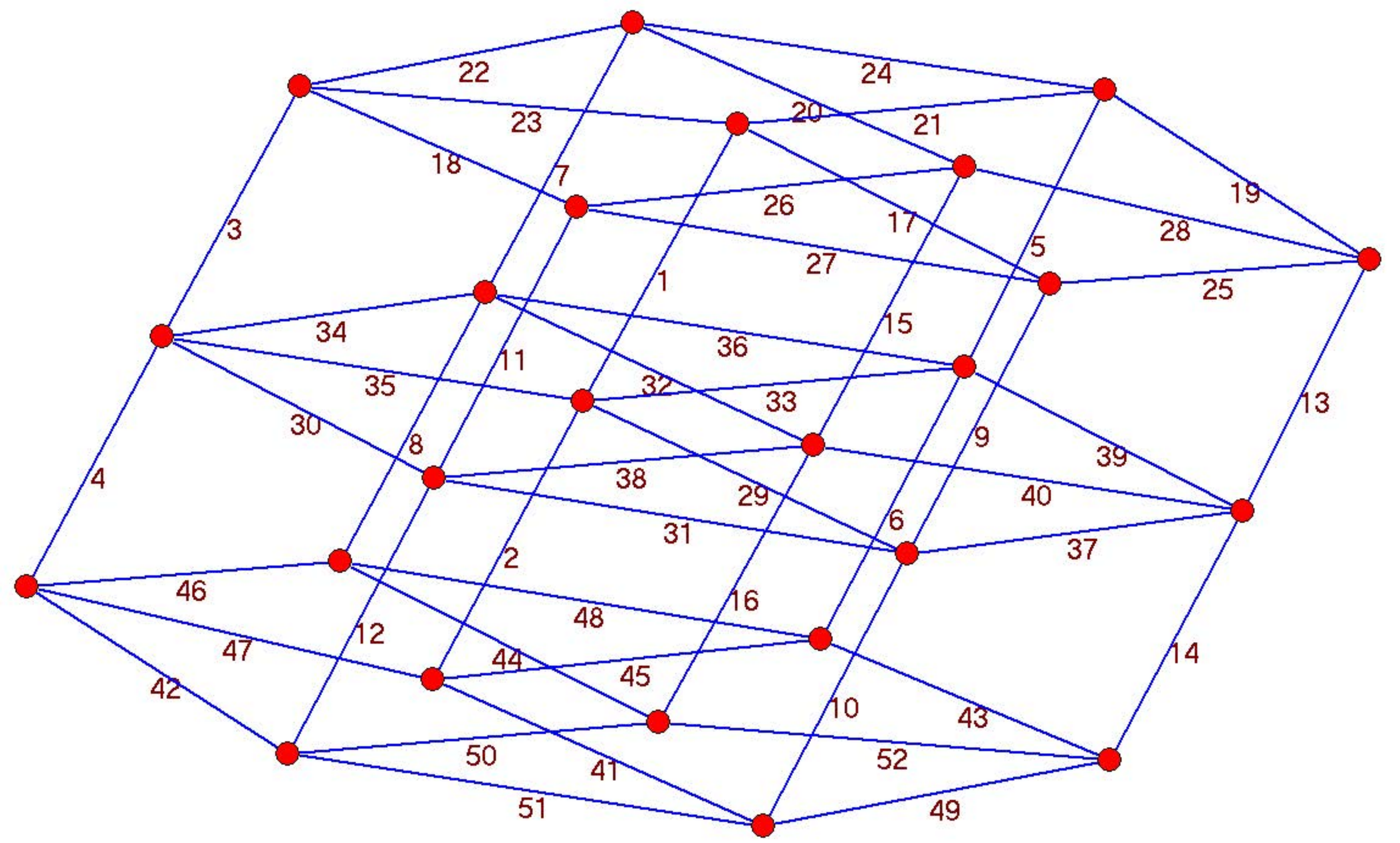}
	 \caption{Vertex-antimagic edge labelling of the cartesian product $P_2^3 \times P_3$.}
	 \label{car24}
\end{center}	
\end{figure}

The algorithm has been implemented in C++ and MATLAB, and the source code can be obtained upon request to the authors. The C++ implementation contains several additional features not described in this paper. For example, it is possible to search for super labellings. The MATLAB implementation also includes functions to convert to$/$from the Simple Graph Format, of MATGRAPH \cite{Schein08}. MATGRAPH was used to generate some of the graphs, and to generate the drawings, in connection with PAJEK \cite{pajek}.
\\
The software described in this paper has provided many useful examples of magic and antimagic graph labellings. For instance, we used our method to check that every wheel $W_k$ for $k$ smaller than 11 admits a vertex-magic total labelling for all feasible values of the magic sum. This completed the solution of vertex-magic total labelling of all wheels, since we know that a wheel $W_n$ has no vertex-magic total labelling if $n > 11$.
\\
We also used our algorithm on trees, and it took a couple of minutes to verify that all trees with less than 10 nodes have an edge-magic super labelling. 
\\
In the case of antimagic labellings, there is a long standing conjecture stating that all connected graphs different from $K_2$ have a vertex-antimagic edge labelling. There is a striking difference between the apparent hardness of this conjecture, and the apparent easiness to find an antimagic labelling for a particular connected graph. In a recent paper \cite{Zhan09}, it was proved that the cartesian product of two path graphs, $P_i$ and $P_j$, with $i$ and $j$ vertices respectively, is antimagic\footnote{Has a vertex-antimagic edge labelling.} if $i > 3$ or $j > 3$. For the case of $i \leq 3$ and $j \leq 3$, no general construction has been found yet. With the aid of our program, we have easily found vertex-antimagic edge labellings for all cartesian products $P_2^r \times P_3^s$, with less than 50 vertices. In Figure \ref{car24} we show an antimagic labelling of the cartesian product $P_2^3 \times P_3$, with 24 vertices and 52 edges. The files for the other graphs are available upon request. 

\begin{figure}[htbp]
\begin{center}
	 	\includegraphics[width=0.5\textwidth]{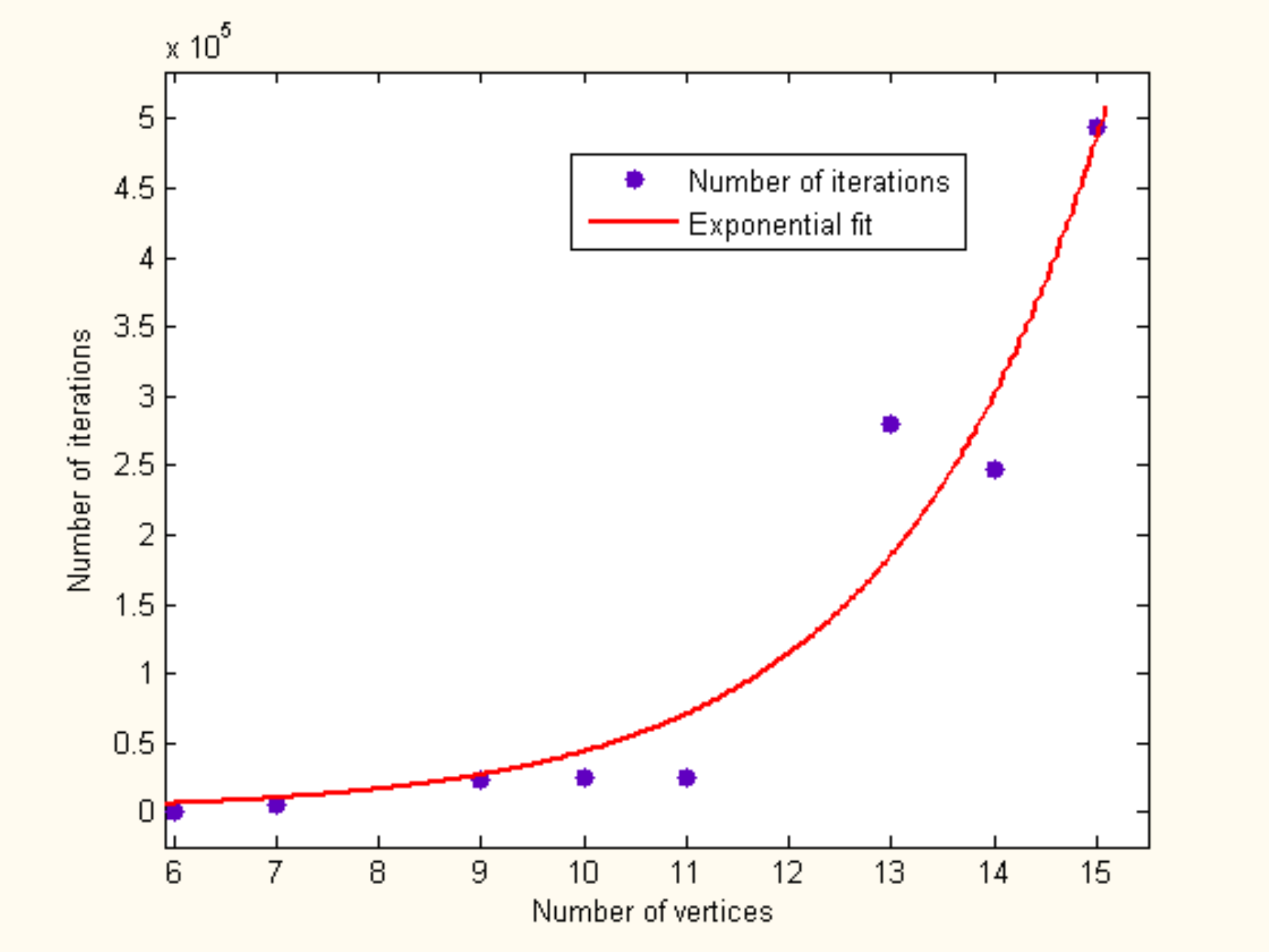}
	 \caption{Number of iterations as a function of $n$.}
	 \label{fit}
\end{center}	
\end{figure}

\begin{figure}[htbp]
\begin{center}
	 	\includegraphics[width=0.5\textwidth]{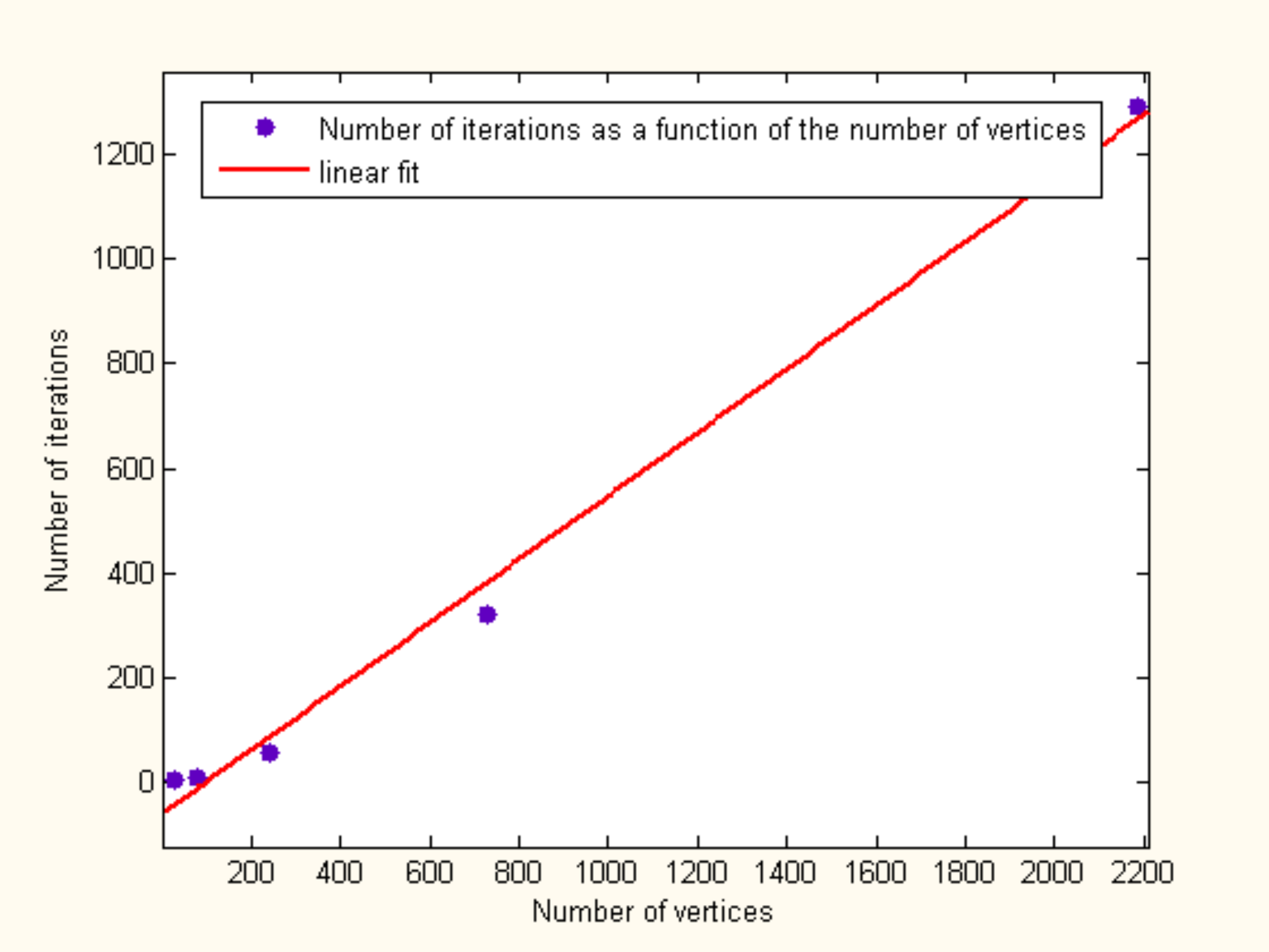}
	 \caption{Number of iterations as a function of $n$.}
	 \label{linearfit}
\end{center}	
\end{figure}

However, the following problem remains open in general: \\
\noindent \textbf{Open Problem}: Is the cartesian product $P_2^r \times P_3^s$ antimagic for all positive integers $r,s$?
\\\\
Evaluating the performance of any simulated annealing algorithm is a difficult task, given the inherent randomness of the process, which leads to great fluctuations in the running time and other parameters of the algorithm. In this particular case we are confronted with additional difficulties: Testing the algorithm on random graphs is out of the question, since the input graph may very well not have a labelling of the desired type, in which case the program will keep running until we decide to stop it. Moreover, there seems to be a great difference in the running time for different types of labellings. For instance, antimagic-type labellings are considerably easier to find than magic ones, in general. To appreciate those variations we have chosen two benchmark problems, one of each type: finding a vertex magic edge superlabelling for the family of complete graphs $K_n$, and finding a vertex-antimagic edge labelling for the family of Cartesian products ${P_3}^k$. 
\\
It is well known that $K_n$ has a vertex magic edge superlabelling for $n \geq 6$, and $n \ne 0 \mbox{ (mod } 4)$. Moreover, the graphs $K_n$ are easy to construct, and there is one for each value of $n$. In principle we could have also chosen the family of graphs $C_n \times C_n$ (Cartesian produc of a cycle with itself), or the family of graphs $Q_d$ (the $d$-dimensional cubes for even $d > 2$), which are also known to have a magic superlabelling, but in these families the number of vertices (and therefore the number of edges) increases too fast, making it impossible to test but for a couple of members of the family. 
\\
We averaged the number of iterations of eight runs for each $6 \leq n \leq 15$, $n \neq 8, 12$, and plotted the results in Figure \ref{fit}. The best fit is the exponential function $y = 357.8 e^{0.4811 x}$, which is consistent with our belief that the problem is NP-hard. 
\\
In the case of vertex-antimagic edge labellings for ${P_3}^k$ we proceeded in the same manner. This problem seems to have a much lower computational complexity, and in fact, constructing the actual graphs took much more time than finding the labellings for them. The difficulty of constructing ${P_3}^k$ prevented us from running any tests for $k > 7$. The results about the number of iterations is shown in Figure \ref{linearfit}. In this case the data can be approximated very well by the linear function $y = 0.653 x - 57.08$. 

\section{Conclusions and future work}
\label{sec:conclusions}

As we have seen, treating the labelling problem as an optimization problem, a very general algorithm can be designed to find a wide range of graph labellings. The objective function of the optimization problem can be plugged into the algorithm according to the particular labelling desired. We believe that this unified approach to treat different types of labelling is our main contribution. Hence, it would be interesting to see if this approach can be extended to other types of labellings, like graceful, sum, totally magic, and $L_{(2,1)}$ labellings. 
\\
We have reviewed different strategies for solving the optimization problem, and we have implemented one of them: simulated annealing. Preliminary computational experiments have shown the algorithm to be effective and relatively efficient, enabling us to find some labellings that were unknown so far. However, more extensive experiments are still required to completely establish the properties of the algorithm. 
\\
In the case of magic labellings, it would be interesting to compare our algorithm with that of \cite{Sun94}. Unfortunately, that algorithm was never implemented, and there is no other program for finding magic labellings, that we know of. A future research task is to implement other algorithms, and devise a set of experiments to compare them all.  
\\
However, it is important to bear in mind that the algorithm in \cite{Sun94} and ours are of a completely different nature, and thus, comparisons between them are to be taken with caution: The algorithm in \cite{Sun94} is an exhaustive search process, that always finds the right answer (if it exists), at the expense of computing time. Our algorithm is a heuristic strategy that attempts to find an answer quickly, but may fail to find it in some cases. Moreover, the algorithm in \cite{Sun94} is only valid for magic labellings, whilst our algorithm is applicable to a broad family of labellings, including supermagic labellings. 
\\
A more meaningful comparison can be done between different metaheuristics for solving the same optimization problem. Hence, it would be interesting to implement another heuristic strategy, such as an evolutionary method, and compare it with the implemented simulated annealing approach. The evolutionary approach looks attractive for various reasons. The danger of every local search technique is that the algorithm can get trapped in a local minimum of the objective function. Simulated annealing can reduce the risk of that outcome, but in some cases it falls short of providing the desired level of robustness. In that case, evolutionary techniques constitute a good alternative. Also, evolutionary algorithms lend themselves very well to parallel implementation, while simulated annealing is essentially a sequential process.


\end{document}